\documentclass[twoside]{article}
\usepackage{amsmath}
\usepackage{amssymb,amsfonts}
\usepackage{graphicx}                 
\newcommand{\bA}{\underline{\underline{A}}}

\newcommand{\bE}{\underline{\underline{E}}}

\newcommand{\bH}{\mathbf{H}}

\newcommand{\bmf}{\mathbf{f}}
\newcommand{\be}{\underline{\underline{e}}}

\newcommand{\ba}{\underline{\underline{a}}}

\newcommand{\bx}{\mathbf{x}}
\newcommand{\by}{\mathbf{y}}

\newcommand{\bu}{\mathbf{u}}
\newcommand{\bv}{\mathbf{v}}

\newcommand{\BE}{\mathbf{E}}
\newcommand{\Bsu}{\boldsymbol{u}}
\newcommand{\Bsv}{\boldsymbol{v}}
\newcommand{\Bsp}{\boldsymbol{p}}
\newcommand{\Bu}{\mathbf{u}}
\newcommand{\Bv}{\mathbf{v}}

\newcommand{\cC}{\mathcal{C}}

\newcommand{\cE}{\mathcal{E}}
\newcommand{\bbR}{\mathbb{R}}

\begin{document}


\title{Isoptic curves of generalized conic sections in the hyperbolic plane}

\author{G\'eza~Csima and Jen\H o~Szirmai, \\
Budapest University of Technology and Economics, \\
Institute of Mathematics,
Department of Geometry Budapest, \\
P.O. Box 91, H-1521 \\
csgeza@math.bme.hu, szirmai@math.bme.hu}
\maketitle
\thanks{}



\begin{abstract}
After having investigated the real conic sections and their isoptic curves in the hyperbolic plane $\bH^2$ we
consider the problem of the isoptic curves of generalized conic sections in the extended hyperbolic plane.

This topic is widely investigated in the Euclidean plane $\BE^2$ (see for example \cite{Lo}), but
in the hyperbolic and elliptic planes there are few results (see \cite{CsSz1}, \cite{CsSz2} and \cite{CsSz3}).
In this paper we recall the former results on isoptic curves in the hyperbolic plane geometry,
and define the notion of the generalized hyperbolic angle between proper and non-proper straight lines, summarize the generalized hyperbolic conic sections classified by K.~Fladt in \cite{KF1} and \cite{KF2} and by E.~Moln\'ar in \cite{M81}.
Furthermore, we determine and visualize the generalized isoptic curves to all hyperbolic conic sections.

We use for the computations the classical
model which are based on the projective interpretation of the hyperbolic geometry and in this manner the isoptic curves
can be visualized on the Euclidean screen of computer.
\end{abstract}

\newtheorem{theorem}{Theorem}[section]
\newtheorem{corollary}[theorem]{Corollary}
\newtheorem{lemma}[theorem]{Lemma}
\newtheorem{exmple}[theorem]{Example}
\newtheorem{defn}[theorem]{Definition}
\newtheorem{rmrk}[theorem]{Remark}
\newenvironment{definition}{\begin{defn}\normalfont}{\end{defn}}
\newenvironment{remark}{\begin{rmrk}\normalfont}{\end{rmrk}}
\newenvironment{example}{\begin{exmple}\normalfont}{\end{exmple}}
\newtheorem{remarque}{Remark}



\section{Introduction}

Let $G$ be one of the constant curvature plane geometries, the Euclidean $\BE^2$, the hyperbolic $\bH^2$, and the elliptic $\cE^2$.
The isoptic curve of a given plane curve $\mathcal{C}$ is the locus of points $P \in G$, where $\mathcal{C}$ is seen under
a given fixed angle $\alpha$ $(0<\alpha <\pi)$.
An isoptic curve formed by the locus of tangents meeting at right angle is called orthoptic curve.
The name isoptic curve was suggested by Taylor in \cite{T}.

In \cite{CMM91} and \cite{CMM96}, the Euclidean isoptic curves of the closed, strictly convex curves are studied, using their support function.
Papers \cite{Kur}, \cite{Wu71-1} and \cite{Wu71-2} deal with Euclidean curves having a
circle or an ellipse for an isoptic curve. Further curves appearing as isoptic curves are well studied in Euclidean plane geometry
$\BE^2$, see e.g. \cite{Lo,Wi}.
Isoptic curves of conic sections have been studied in \cite{H} and \cite{S}. There are results for Bezier curves as well, see \cite{Kunk}.
A lot of papers concentrate on the properties of the isoptics, e.g. \cite{M,MM,MM2}, and the references given there. There are some generalization of the isoptics as well \textit{e.g.} equioptic curves in \cite{O} or secantopics in \cite{Skrzy}

In the case of hyperbolic plane geometry there are only few results.
The isoptic curves of the hyperbolic line segment and proper conic sections are determined by the authors in \cite{CsSz1}, \cite{CsSz2} and \cite{CsSz3}.

The isoptics of conic sections in elliptic geometry $\cE^2$ are determined by the authors in \cite{CsSz3}.

In the papers \cite{KF1} and \cite{KF2} K.~Fladt determined the equations of the generalized conic sections in the hyperbolic plane
using algebraic methods and in \cite{M81} E.~Moln\'ar classified them with synthetic methods.

Our goal in this paper is to generalize our method described in \cite{CsSz3}, that is based on the projective interpretation of hyperbolic plane geometry,
to determine the isoptic curves of the generalized hyperbolic conics and visualize them for some angles. Therefore we study and recall
the notion of the angle between proper and non-proper straight lines using the results of the papers \cite{BH}, \cite{GH} and \cite{V}.

\section{The projective model}

For the $2$-dimensional hyperbolic plane $\bH^2$ we use the projective model in Lorentz space
$\BE^{2,1}$ of signature $(2,1)$, i.e.~$\BE^{2,1}$ is
the real vector space $\mathbf{V}^{3}$ equipped with the bilinear
form of signature $(2,1)$
\begin{equation}
\langle ~ \mathbf{x}, ~ \mathbf{y} \rangle = x^1y^1+ x^2 y^2-x^3 y^3
\end{equation}
where the non-zero vectors
$\mathbf{x} = (x^1 , x^2 , x^3)^T $ and $\by= (y^1 , y^2 , y^3)^T \in \mathbf{V}^{3},$
are determined up to real factors and they represent points $X = \mathbf{x}\mathbb{R}$ and $Y = \mathbf{y}\mathbb{R}$ of $\bH^2$ in
$\mathbb{P}^2(\mathbb{R})$.
The proper points of $\bH^2$ are represented as the interior of the absolute conic
\begin{equation}
AC=\{\mathbf{x}\mathbb{R} \in\mathcal{P}^2 | \langle ~ \mathbf{x},~\mathbf{x} \rangle =0 \}=\partial \bH^2
\end{equation}
in real projective space $\mathbb{P}^2(\mathbf{V}^{3}, \mbox{\boldmath$V$}\!_{3})$.
All proper interior point $X \in \bH^2$ are characterized by
$\langle ~ \mathbf{x},~\mathbf{x} \rangle < 0$. The points on the boundary $\partial \bH^2 $ in
$\mathcal{P}^2$ represent the absolute points at infinity of $\bH^2 $.
Points $Y$ with $\langle ~ \mathbf{y},~\mathbf{y} \rangle >0$ are called outer or non-proper points
of $\bH^2 $.

The point $Y=\mathbf{y}\mathbb{R}$ is said to be conjugate to $X=\mathbf{x}\mathbb{R}$ relative to $AC$ when $\langle ~\mathbf{x},~\mathbf{y} \rangle =0$.

The set of all points conjugate
to $X=\mathbf{x}\mathbb{R}$ forms a projective (polar) line
\begin{equation}
pol(X):=\{\mathbf{y}\mathbb{R} \in \mathbb{P}^2 | \langle ~ \mathbf{x},~\mathbf{y} \rangle =0 \}.
\end{equation}

Hence the bilinear form to $(AC)$ by (1) induces a bijection
(linear polarity $\mathbf{V}^{3} \rightarrow
\mbox{\boldmath$V$}\!_{3})$
from the points of $\mathcal{P}^2$
onto its lines (hyperplanes in general).

Point $X=\mathbf{x}\mathbb{R}$ and the hyperplane $u=\mathbb{R}\mbox{\boldmath$u$}$ are called incident if the value of the
linear form $\mbox{\boldmath$u$}$ on the vector $\mathbf{x}$ is equal
to zero; i.e., $\mbox{\boldmath$u$}\mathbf{x}=0$ ($\mathbf{x} \in \
\mathbf{V}^{3} \setminus \{\mathbf{0}\}, \ \mbox{\boldmath$u$} \in
\mbox{\boldmath$V$}_{3} \setminus \{\mbox{$\boldsymbol{0}$}\}$).
In this paper we set the sectional curvature of $\bH^2$,
$K=-k^2$, to be $k=1$.

The distance $d(X,Y)$ of two proper points
$X=\mathbf{x}\mathbb{R} $ and $Y=\mathbf{y}\mathbb{R}$ can be calculated with appropriate representant vectors by the formula (see e.g. \cite{MSz}) :
\begin{equation}
\cosh{{d(X,Y)}}=\frac{-\langle ~ \mathbf{x},~\mathbf{y} \rangle }{\sqrt{\langle ~ \mathbf{x},~\mathbf{x} \rangle
\langle ~ \mathbf{y},~\mathbf{y} \rangle }}.
\end{equation}

For the further calculations, let denote $\bu$ the pole of the straight line $u=\mathbb{R}\mbox{\boldmath$u$}$.
It is easy to prove, that if $\mbox{\boldmath$u$}=(u_1,u_2,u_3)$ then $\bu=(u_1,u_2,-u_3)$.
And follows that if $u=\bbR\mbox{\boldmath$u$}$ and $v=\bbR\mbox{\boldmath$v$}$
then $\left\langle \bu,\bv\right\rangle=\left\langle \mbox{\boldmath$u$},\mbox{\boldmath$v$} \right\rangle$.

\subsection{Generalized angle of straight lines}

Having regard to the fact that the majority of the generalized conic sections have ideal and outer tangents as well, it is inevitable to
introduce the generalized concept of the hyperbolic angle.
In the extended hyperbolic plane there are three classes of lines by the number of common points with the absolute conic $AC$ (see (2)):
\begin{enumerate}
\item The straight line $u=\mathbb{R}\mbox{\boldmath$u$}$ is \textit{proper} if card${(u\cap AC)}=2$ $\Leftrightarrow $
$\langle \mbox{\boldmath$u$},~\mbox{\boldmath$u$} \rangle>0$.
\item The straight line $u=\mathbb{R}\mbox{\boldmath$u$}$ is \textit{non-proper} if card${(u\cap AC)} < 2.$
\begin{enumerate}
\item If card${(u\cap AC)}=1$ $\Leftrightarrow $ $\langle \mbox{\boldmath$u$},~\mbox{\boldmath$u$} \rangle=0$  then $u=\mathbb{R}\mbox{\boldmath$u$}$
is called \textit{boundary} straight line.
\item If card${(u\cap AC)} = 0$ $\Leftrightarrow $ $\langle \mbox{\boldmath$u$},~\mbox{\boldmath$u$} \rangle<0$ then
$u=\mathbb{R}\mbox{\boldmath$u$}$ is called \textit{outer} straight line.
\end{enumerate}
\end{enumerate}
We define the generalized angle between straight lines using the results of the papers \cite{BH}, \cite{GH} and \cite{V} in the projective model.
\begin{definition}
\begin{enumerate}

\item Suppose that $u=\mathbb{R}\mbox{\boldmath$u$}$ and $v=\mathbb{R}\mbox{\boldmath$v$}$ are both \textit{proper} lines. 
\begin{enumerate}
\item If $\langle \mbox{\boldmath$u$},\mbox{\boldmath$u$} \rangle \langle \mbox{\boldmath$v$},\mbox{\boldmath$v$} \rangle -
\langle \mbox{\boldmath$u$},\mbox{\boldmath$v$} \rangle^2 >0$ then they intersect in a proper point and their angle
$\alpha (\mbox{\boldmath$u$},\mbox{\boldmath$v$})$ can be measured by
				\begin{equation}
				 \cos{\alpha}=\frac{\pm\langle \Bsu,\Bsv \rangle}{\sqrt{
				 \langle \Bsu,\Bsu\rangle \langle \Bsv,\Bsv \rangle}}.
				 \label{two}
				\end{equation}
\item If $\langle \mbox{\boldmath$u$},\mbox{\boldmath$u$} \rangle \langle \mbox{\boldmath$v$},\mbox{\boldmath$v$} \rangle -
\langle \mbox{\boldmath$u$},\mbox{\boldmath$v$} \rangle^2 <0$ then they intersect in a non-proper point and their angle is the      				\textit{length of their normal transverse} and it can be calculated using the formula below:
\begin{equation}
				 \cosh{\alpha}=\frac{\pm\langle \Bu,\Bv \rangle}{\sqrt{
				 \langle \Bu,\Bu\rangle \langle \Bv,\Bv \rangle}}.
				 \label{three}
				\end{equation}
	\item If $\langle \mbox{\boldmath$u$},\mbox{\boldmath$u$} \rangle \langle \mbox{\boldmath$v$},\mbox{\boldmath$v$} \rangle -
\langle \mbox{\boldmath$u$},\mbox{\boldmath$v$} \rangle^2 =0$  then they intersect in a boundary point and their angle is $0$.			
 \end{enumerate}
\item Suppose that $u= \mathbb{R} \Bsu $ and $v= \mathbb{R} \Bsv $ are both \textit{outer} lines of $\bH^2$. The angle of these lines will be the
\textit{distance of their poles} using the formula (6).
\item Suppose that $u= \mathbb{R} \Bsu $ is a \textit{proper} and $v= \mathbb{R} \Bsv $ is an \textit{outer} line.
Their angle is defined as the \textit{distance of the pole of the outer line to the real line} and can be computed by
\begin{equation}\label{four}
\sinh\alpha=\frac{\pm\langle \Bu,\Bv \rangle}{\sqrt{-
\langle \Bu,\Bu\rangle \langle \Bv,\Bv \rangle}}.
\end{equation}
\item Suppose that at least one of the straight lines $u= \mathbb{R} \Bsu $ and $v= \mathbb{R} \Bsv $ is \textit{boundary}
line of $\bH^2$. If the other line fits the boundary point, the angle cannot be defined, otherwise it is infinite.
\end{enumerate}
\end{definition}
\begin{remark}
In the previous definition we fixed that except case 1. (a) we use real \textit{distance type} values instead of \textit{complex angles} which arise in other cases. The $\pm$ on the right sides are justifiable because we consider complementary angles \textit{i.e.} $\alpha$ and $\pi-\alpha$ together.

\label{rmr1}
\end{remark}
\section{Classification of generalized conic sections on the hyperbolic plane in dual pairs}


In this section we will summarize and extend the results of K.~Fladt (see \cite{KF1} and \cite{KF2}) about the generalized conic sections on the extended hyperbolic plane.

Let us denote a point with $\bx$ and a line with $\Bsu$. Then the absolute conic (AC) can be defined as a point conic with the $\bx^T\be\bx=0$ quadratic form where $\be=\mathrm{diag}\left\{1,1,-1\right\}$ or due to the absolute polarity as line conic with $\Bsu \bE\Bsu^T=0	$ where $\bE=\be^{-1}=\mathrm{diag}\left\{1,1,-1\right\}$.

Similarly to the Euclidean geometry we use the well-known quadratic form
$$\mathbf{x}^T \ba\mathbf{x}=a_{11} x^1x^1+ a_{22} x^2x^2 + a_{33} x^3x^3+2 a_{23} x^2x^3+2a_{13} x^1x^3+2 a_{12} x^1x^2=0 $$
where $\det{a}\neq0$ for a non-degenerate point conic and
$$\Bsu \bA\Bsu^T=A^{11} u_1u_1+ A^{22} u_2u_2 + A^{33} u_3u_3+2 A^{23} u_2u_3+2A^{13} u_1u_3+2 A^{12} u_1u_2=0 $$
where $\bA=\ba^{-1}$ for the corresponding line conic defined by the tangent lines of 
 the previous point conic. Using the polarity $\bx=\bA\Bsu^T$ and $\Bsu^T=\ba\bx$ follow since $\Bsu\bx=0$.

Consider a one parameter conic family of our point conic with the (AC), defined by
\begin{equation}
\bx^T(\ba+\rho \be)\bx=0.
\notag
\end{equation}

Since the characteristic equation $\Delta(\rho):=\det(\ba+\rho \be)$ is an odd degree polynomial, this conic pencil has at least one real degenerate element $(\rho_1)$, which consists of at most two point sequences with holding lines $\Bsp_1^1$ and $\Bsp_1^2$ called asymptotes. Therefore we get a product
\begin{equation}
\bx^T(\ba+\rho_1 \be)\bx=(\Bsp_1^1\bx)^T(\Bsp_1^2\bx)=\bx^T((\Bsp_1^1)^T\Bsp_1^2)\bx=0
\notag
\end{equation}
with occasional complex coordinates of the asymptotes. Each of these two asymptotes has at most two common points with the (AC) and with the conic as well. Thus, the at most 4 common points with at most 3 pairs of asymptotes can be determined through complex coordinates and elements according to the at most 3 different eigenvalues $\rho_1$, $\rho_2$ and $\rho_3$.

In complete analogy with the previous discussion in dual formulation we get that the one parameter conic family of a line conic with (AC) has at least one degenerate element $(\sigma^1)$ which contains two line pencils at most with occasionally complex holding points $\bmf_1^1$ and $\bmf_2^1$ called foci.

\begin{equation}
\Bsu(\bA+\sigma^1 \bE)\Bsu^T=(\Bsu \bmf_1^1)(\Bsu \bmf_2^1)^T=\Bsu(\bmf_1^1(\bmf_2^1)^T)\Bsu^T=0
\notag
\end{equation}

For each focus at most two common tangent line can be drawn to $AC$ and to our line conic. Therefore, at most four common tangent lines with at most three pairs of foci can be determined maybe with complex coordinates to the corresponding eigenvalues $\sigma^1$, $\sigma^2$ and $\sigma^3$.

Combining the previous discussions with \cite{KF1} and \cite{M81} the classification of the conics on the extended hyperbolic plane can be obtained in dual pairs.

First, our goal is to find an appropriate transformation, so that the resulted normalform characterizes the conic \textit{e. g.} the straight line $x^1=0$ is a symmetry axis of the conic section  ($a_{31}=a_{12}=0$).
Therefore we take a rotation around the origin $O(0,0,1)^T$ and a translation parallel with $x^2=0$.

As it used before, the characteristic equation$$
\Delta(\rho)=\det(\ba+\rho\be)=
\det\left(\begin{array}{ccc}
	a_{11}+\rho & a_{12} & a_{13}\\
	a_{21} & a_{22}+\rho & a_{23}\\
	a_{31} & a_{32} & a_{33}-\rho
\end{array}\right)=0$$ has at least one real root denoted by $\rho_1$.

This is helpful to determine the exact transformation if  the equalities $\rho_1=\rho_2=\rho_3$  not hold.
That case will be covered later. With this transformations we obtain the normalform
\begin{equation}
\rho_1 x^1x^1+a_{22}x^2x^2+2a_{23}x^2x^3+a_{33}x^3x^3=0.
\label{nf}
\end{equation}
\textit{In the following we distinguish $3$ different cases according to the other two roots:}
%

\begin{enumerate}
\item[1.] \textit{Two different real roots} \\ 
Then the monom $x^2x^3$ can be eliminated from the equation above, by translating the conic parallel with $x^1=0$. The final form of the conic equation in this case, called \textbf{central conic} section:

\[\rho_1 x^1x^1+\rho_2 x^2x^2-\rho_3 x^3x^3=0.\]

Because our conic is non-degenerate $\rho_3x^3\neq0$ follows and with the notations $a=\frac{\rho_1}{\rho_3}$ and $b=\frac{\rho_2}{\rho_3}$ our matrix can be transformed into $\ba=diag\left\{a,b,-1\right\}$, where $a\leq b$ can be assumed. The equation of the dual conic can be obtained using the polarity $\bE$ respected to (AC) by $\bE~\bA~\bE^{-1}=diag\{\frac{1}{a},\frac{1}{b},-1\}$.
By the above considerations we can give an overview of the generalized central conics with representants:
\begin{theorem}
If the conic section has the normalform $a x^2+b y^2=1$ then we get the following types of central conic sections (see Figures 1-3):
\begin{enumerate}
	\item \textnormal{Absolute conic:} \ \ \ \ \ \ \ \ \ \ \ \ \ \ \ \ \ \ \ \ \ \ \ \ \ \ \ \ \ \ \ \ \ \ \ \ \ \ \ \ \ $a=b=1$
	\item \begin{enumerate}
				\item \textnormal{Circle:} \ \ \ \ \ \ \ \ \ \ \ \ \ \ \ \ \ \ \ \ \ \ \ \ \ \ \ \ \ \ \ \ \ \ \ \ \ \ \ \ \ \ \ \ \ \ \ $1<a=b$
				\item \textnormal{Circle enclosing the absolute:} \ \ \ \ \ \ \ \ \ \ \ \ \ \ \ \ \ \ \  $a=b<1$
				\end{enumerate}
	\item \begin{enumerate}
				\item \textnormal{Hypercycle:} \ \ \ \ \ \ \ \ \ \ \ \ \ \ \ \ \ \ \ \ \ \ \ \ \ \ \ \ \ \ \ \  \ \ \ \ \ \ \ \, $1=a<b$
				\item \textnormal{Hypercycle enclosing the absolute: }\ \ \ \ \ \ \ \ \ \ $0<a<1=b$
				\end{enumerate}
	\item \textnormal{Hypercycle excluding the absolute:}\ \ \ \ \ \ \ \ \ \ \ \ \ \ \ \,$a<0<1=b$
	\item \textnormal{Concave hyperbola:}  \ \ \   \ \ \ \ \ \ \ \ \ \ \ \ \ \ \ \ \ \ \, \ \ \ \ \ \ \ \ \ \,  $0<a<1<b$
	\item \begin{enumerate}
				\item \textnormal{Convex hyperbola:}   \ \ \     \ \ \ \ \ \ \ \ \ \ \ \ \ \ \ \ \ \ \ \ \ \ \ \ \ \    $a<0<1<b$
				\item \textnormal{Hyperbola excluding the absolute:}\ \ \ \ \ \ \ \ \ \ \, $a<0<b<1$
				\end{enumerate}
	\item \begin{enumerate}
			  \item \textnormal{Ellipse:}	 \ \ \ \ \ \ \ \ \ \ \ \ \ \ \ \ \ \ \ \ \ \ \ \ \ \ \ \ \ \ \ \ \ \ \ \: \ \ \ \ \ \ \ \, $1<a<b$
				\item \textnormal{Ellipse enclosing the absolute:}  \ \ \: \ \ \ \ \ \ \ \ \ \ \,  $0<a<b<1$
				\end{enumerate}
	\item \textnormal{empty:}\ \ \ \ \ \ \ \ \ \; \ \ \ \ \ \ \ \ \ \ \ \ \ \ \ \ \ \ \ \ \ \ \ \ \ \ \ \ \ \ \ \: \ \ \ \ \ \ \ \, $a\leq b\leq0$
\end{enumerate}
where either the conic and its dual pair lies in the same class or (i) and (ii) are dual pairs with $a'=\frac{1}{a}$ and $b'=\frac{1}{b}$.
\end{theorem}
\item[2.] \textit{Coinciding real roots}\\
The last translation cannot be enforced but it can be proved that $\rho_2=\rho_3=\frac{a_{33}-a_{22}}{2}$ follows.
With some simplifications of the formulas in \cite{KF1} we obtain the normalform of the so-called generalized {\bf parabolas}.
\begin{theorem}
The parabolas have the normalform $ax^2+(b+1)y^2-2y=b-1$ and the following cases arise (see Figure 4-6):
\begin{enumerate}
	\item \begin{enumerate}
				\item \textnormal{Horocycle:} \ \ \ \ \ \ \ \ \ \ \ \ \ \ \ \ \ \ \ \ \ \ \ \ \ \ \ \ \ \ \ \ \ \ \ \ \ \ \ \ \ \ \ \ \ \ \ $0<a=b$
				\item \textnormal{Horocycle enclosing the absolute:} \ \ \ \ \ \ \ \ \ \ \ \ \ \ \ \ \ \ \  $a=b<0$
				\end{enumerate}
	\item \begin{enumerate}
				\item \textnormal{Elliptic parabola:} \ \ \ \ \ \ \ \ \ \ \ \ \ \ \ \ \ \ \ \ \ \ \ \ \ \ \ \ \ \ \ \ \ \ \ \ \ \ $0<b<a$
				\item \textnormal{Parabola enclosing the absolute:} \ \ \ \ \ \ \ \ \ \ \ \ \ \ \ \ \ \ \ \ $b<a<0$
				\end{enumerate}
	\item \begin{enumerate}
				\item \textnormal{Two sided parabola:} \ \ \ \ \ \ \ \ \ \ \ \ \ \ \ \ \ \ \ \ \ \ \ \ \ \ \ \ \  \ \ \ \ \ \ $a<b<0$
				\item \textnormal{Concave hyperbolic parabola:} \ \ \ \ \ \ \ \ \ \ \ \ \ \ \ \ \ \ \ \ \ \ \ \  $0<a<b$
				\end{enumerate}
	\item \begin{enumerate}
				\item \textnormal{Convex hyperbolic parabola:}  \ \ \ \ \ \ \ \ \ \ \ \ \ \ \ \ \ \ \ \ \ \ \ \ \ $a<0<b$
				\item \textnormal{Parabola excluding the absolute: } \ \ \ \ \ \ \ \ \ \ \ \ \ \ \ \ \ \ \ $b<0<a$
				\end{enumerate}
\end{enumerate}
where all (i) and (ii) are dual pairs with parameters $a'=-\frac{b^2}{a}$ and $b'=-b$.
\end{theorem}
\item[3.] \textit{Two conjugate complex roots}\\
Then the last translation cannot be performed to eliminate the monom $x^2x^3$ but we can eliminate the monom $x^3x^3$ by an appropriate transformation described in \cite{KF1}. Shifting to inhomogeneous coordinates and simplifying the coefficients we obtain:
\begin{theorem}
The so-called \textbf{semi-hyperbola} has the normalform $a x^2+2 b y^2-2y=0$ where $\left|b\right|<1$ and its dual pair is projectively equivalent with another semi-hyperbola with $a'=\frac{1}{a}$ and $b'=-b$ (see Figure 7).
\end{theorem}

\item[4.]

Overviewing the above cases only one remains, when the conic has no symmetry axis at all and $\rho_1=\rho_2=\rho_3$.
Ignoring further explanations we claim the following theorem:
\begin{theorem}
If the conic has the normalform
$(1-x^2-y^2)+2a y(x+1)=0$ where $a>0$ then it is called \textbf{osculating parabola}. Its dual is also an osculating parabola by a convenient reflection (see Figure 8)
.
\end{theorem}
\end{enumerate}

\section{Isoptic curves of generalized hyperbolic conics}
\subsection{Isoptic curves of central conics}

In this section we will extend the algorithm for determining the isoptic curves of conic sections described in \cite{CsSz3} for
generalized hyperbolic conic sections. We have to apply the generalized notion of angle, therefore the equation of the isoptic curve
will not be given by a implicit formula but the isoptic curve consist of some piecewise continuous arcs that will be given by implicit equations.

First we determine the equations of the tangent lines through a given external point $P$ to a given conic section $\cC$. We will use not only the point conic but the corresponding line conic as well. The algorithm do not need the points of tangency but of course they can be determined from the tangents. This strategy provides a general simplification in the other cases as well without all details in this work.


Let the external point $P$ be given with homogeneous coordinates $(x,y,1)^T$. 
If a central conic is given by $\ba$, then the corresponding line conic is defined by $\bA$, where $\ba^{-1}=\bA$. Now, we know that $P$ fits on the tangent lines $u= \mathbb{R} \Bsu $ and $v= \mathbb{R} \Bsv $ where $\Bsu=(u_1,u_2,1)$ and $\Bsv=(v_1,v_2,1)$ furthermore  $u$ and $v$ satisfy the equation of the line conic.
$$
\begin{array}{lr}
\left.
\begin{array}{lcl}     
    u_1 x+u_2 y+1&= &0\\
   \frac{u_1^2}{a}+\frac{u_2^2}{b}-1&=&0
\end{array}\right\}
&
\left.
\begin{array}{lcl}     
   v_1 x+v_2 y+1&= &0\\
   \frac{v_1^2}{a}+\frac{v_2^2}{b}-1&=&0
\end{array}\right\}
\end{array}
$$


Solving the above systems we obtain the coordinates of the straight lines $u$ and $v$:
\begin{equation}
\begin{array}{ll}
\begin{array}{l}
   u_1=-\frac{a x+\sqrt{a b y^2 \left(a x^2+b y^2-1\right)}}{a x^2+b y^2}\\
   u_2=\frac{-b y^2+x \sqrt{a b y^2 \left(a x^2+b y^2-1\right)}}{a x^2 y+b y^3}
\end{array}
&
\begin{array}{l}
   v_1=\frac{-a x+\sqrt{a b y^2 \left(a x^2+b y^2-1\right)}}{a x^2+b y^2}\\
   v_2=-\frac{b y^2+x \sqrt{a b y^2 \left(a x^2+b y^2-1\right)}}{a x^2 y+b y^3}
\end{array}
\end{array}
\end{equation}

Of course, $a b y^2 \left(a x^2+b y^2-1\right)\geq0$ must hold otherwise $P(x,y,1)^T$ is not an external point.


We get the exact formula of the more parted isoptic curve related to the central conic sections (see Theorem 3.1) using
the definition of the generalized angle (see Definition 2.1). The {\it compound isoptic} can be given by classifying
the straight lines $u= \mathbb{R} \Bsu $ and $v= \mathbb{R} \Bsv $ according to their poles. We summarize our result in the following:

\begin{theorem}
Let a central conic section be given by its equation $a x^2+b y^2=1$ (see Theorem 3.1).
Then the compound $\alpha$-isoptic curve $(0<\alpha<\pi)$ of the considered conic has the equation
$$
%
%
%
\frac{\left(a \left((b+1) x^2-1\right)+(a+1) b y^2-b\right)^2}{\left| (a-1)^2 b^2 y^4+2 (a-1) b \left(b+a \left((b-1) x^2-1\right)\right) y^2+\left(a
   (b-1) x^2+a-b\right)^2\right|}=
$$
\newline
$$
=\left\{\begin{array}{r}
\begin{array}{ll}
 \cosh^2 (\alpha ), & a b y^2 \left(a x^2+b y^2-1\right)\geq 0~\land \\
 & \left(1-u_1^2-u_2^2\right) \left(1-v_1^2-v_2^2\right)>0~\land\\
 & x^2+y^2>1 \\
 & \\
 \cos^2 (\alpha ), & a b y^2 \left(a x^2+b y^2-1\right)\geq 0~\land\\
 &x^2+y^2<1\\
 & \\
 \sinh^2 (\alpha ), & a b y^2 \left(a x^2+b y^2-1\right)\geq 0~\land \\
 & \left(1-u_1^2-u_2^2\right) \left(1-v_1^2-v_2^2\right)<0,
\end{array}
\end{array}\right.
$$
wherein $u_{1,2}$ and $v_{1,2}$ are derived by \textnormal{(8)}.
\end{theorem}

In Figures 1-3 we visualize the isoptic curves of central conic sections. The foremost figure shows, how different types of isoptics arise due to common tangents. We indicated the Cayley-Klein model circle with black, the conic with dashed line and the isoptic is shaded.

\begin{figure}[p]
\centering
\includegraphics[scale=0.65]{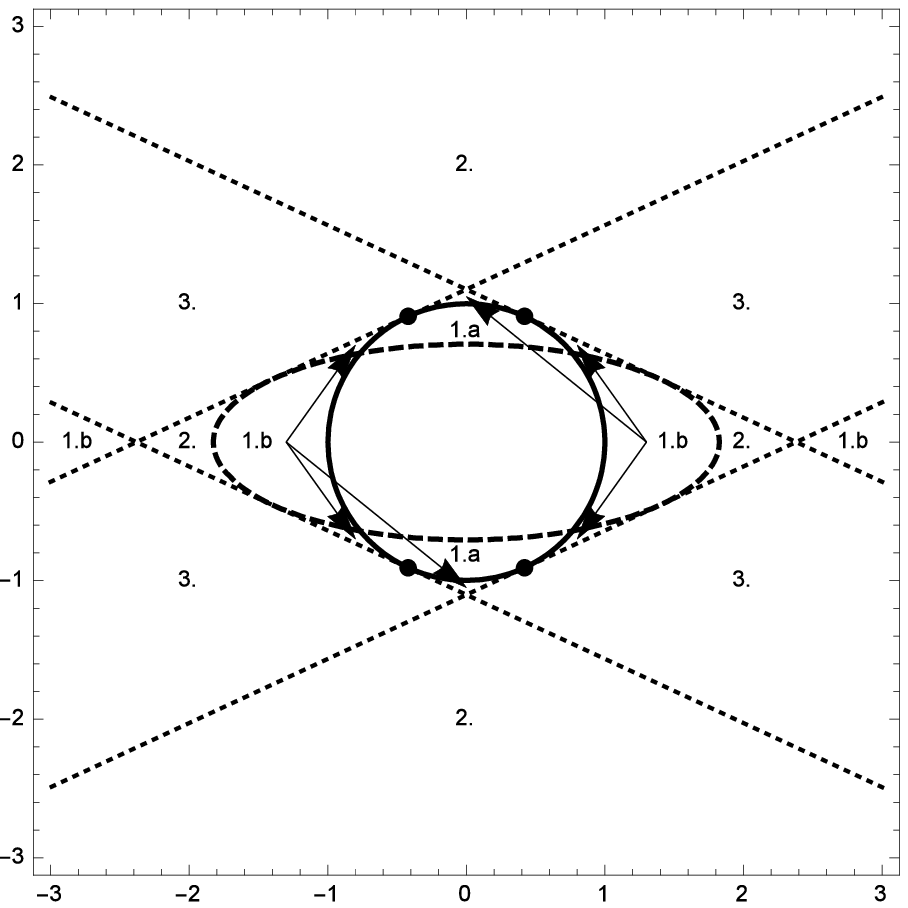}\includegraphics[scale=0.65]{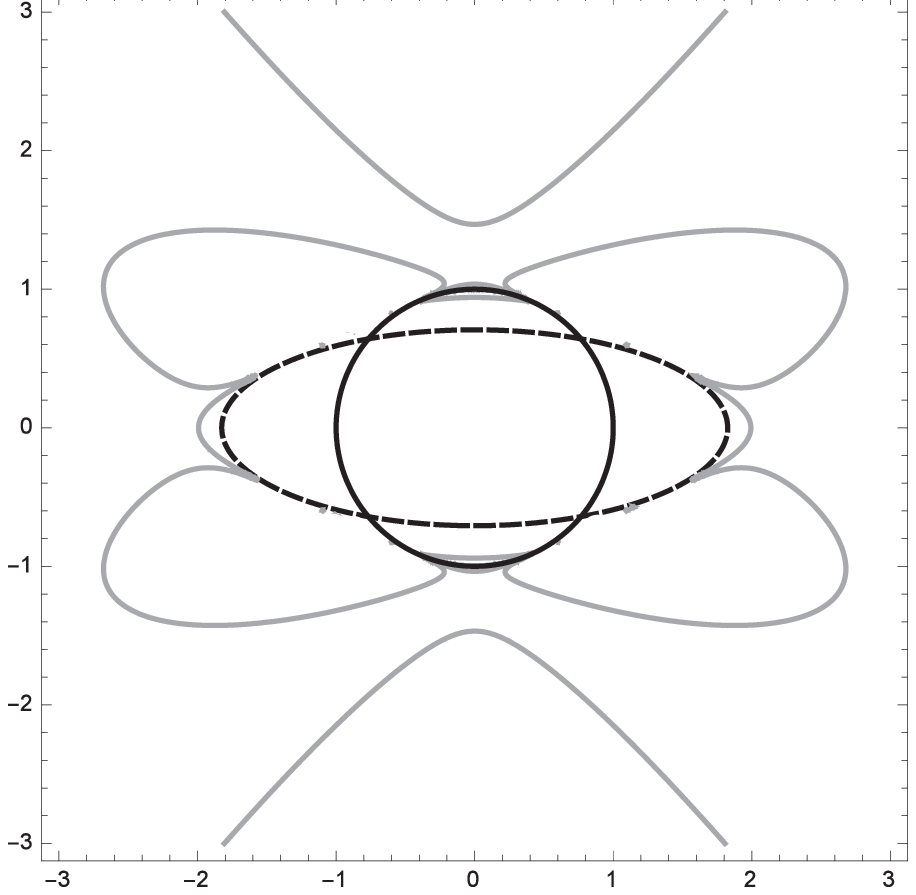}
\caption{Domains for concave hyperbola with notations of Definition 2.1 (left) \newline Concave hyperbola(right): $a=0.3$, $b=2$, $\alpha=\frac{\pi}{2}$}
%
\centering
\includegraphics[scale=0.65]{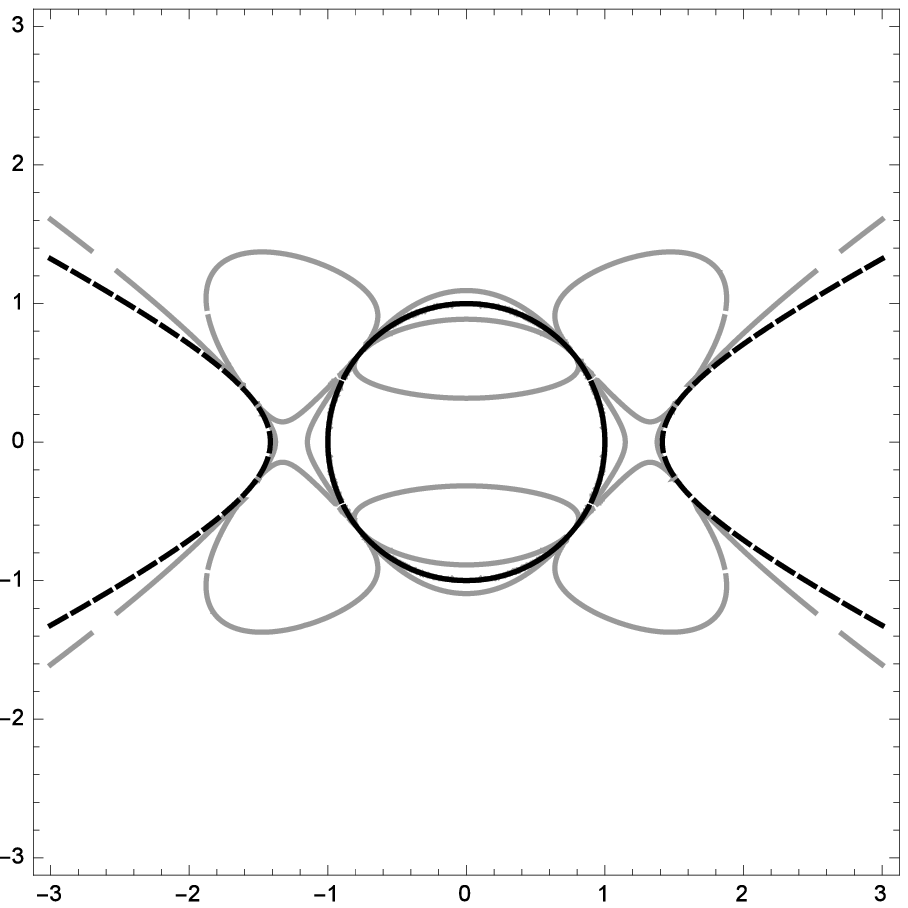}\includegraphics[scale=0.65]{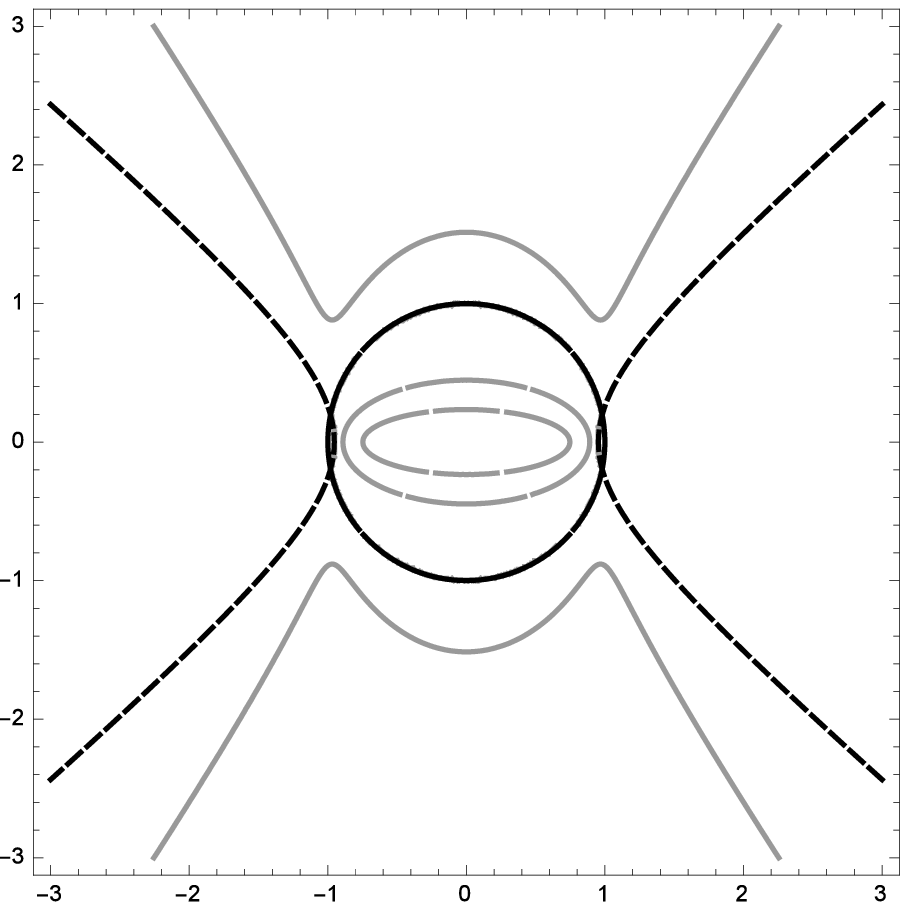}
\caption{Hyperbola excluding the absolute (left): $a=0.5$, $b=-2$, $\alpha=\frac{\pi}{3}$\newline Convex hyperbola(right): $a=1.1$, $b=-1.5$, $\alpha=\frac{19\pi}{36}$}
%
\centering
\includegraphics[scale=0.65]{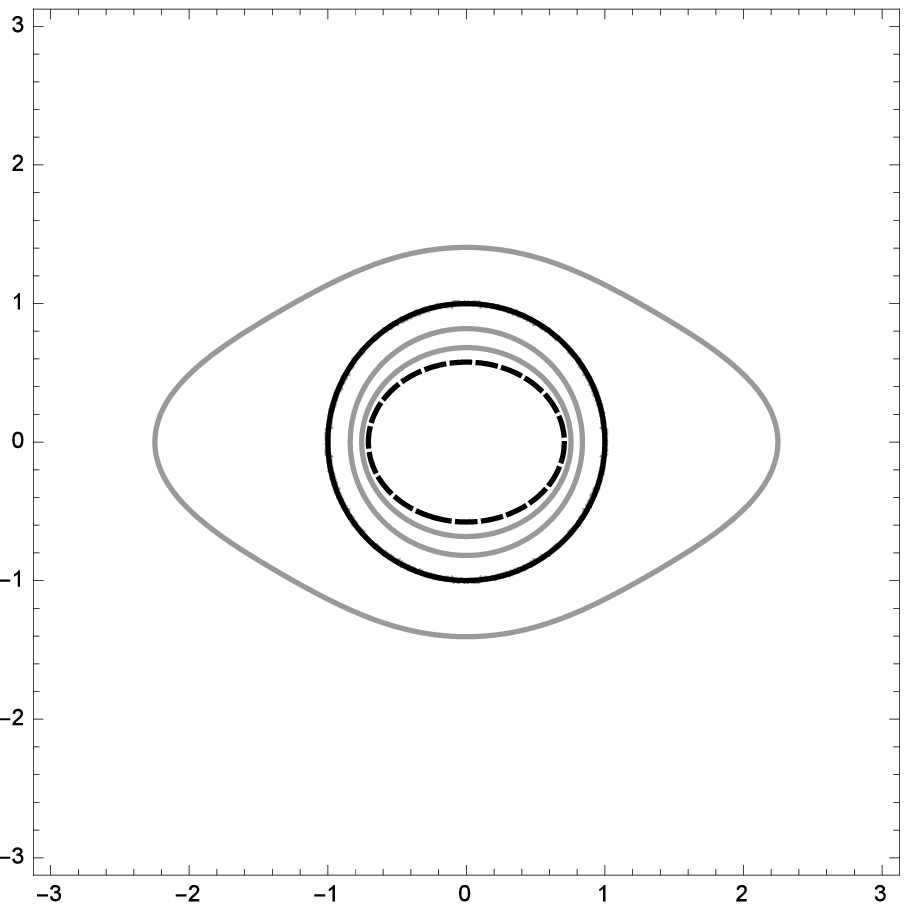}\includegraphics[scale=0.65]{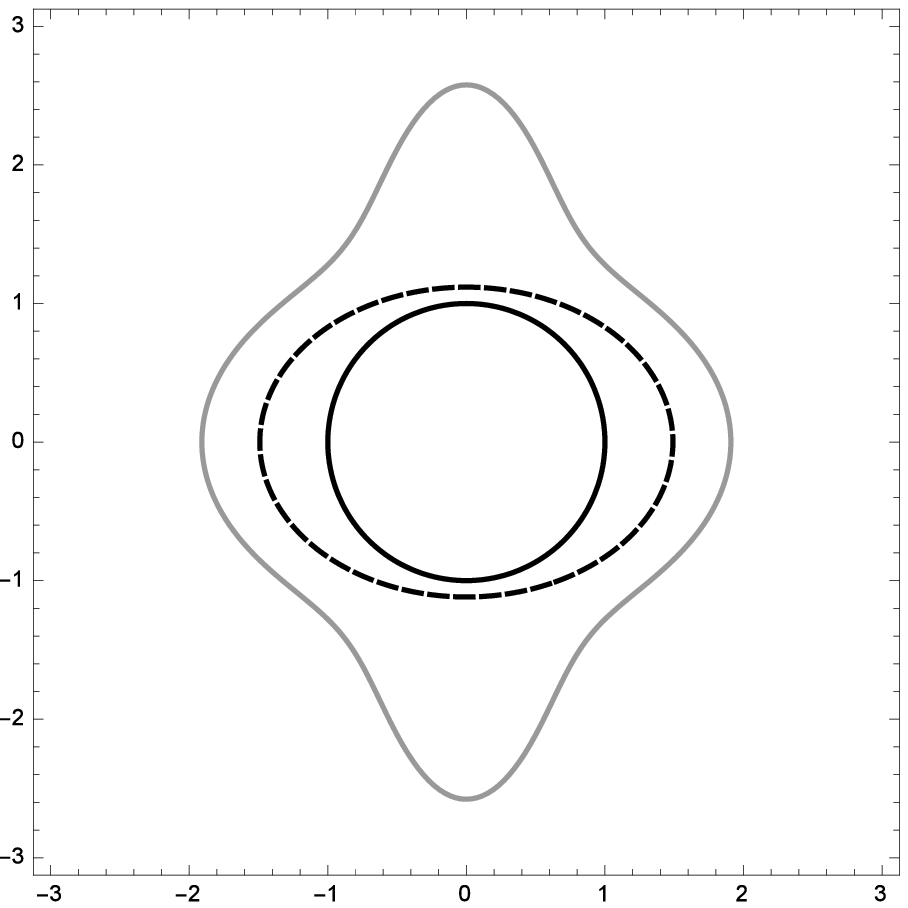}
\caption{Ellipse(left): $a=2$, $b=3$, $\alpha=\frac{7\pi}{18}$,\newline Ellipse enclosing the absolute (right): $a=0.45$, $b=0.8$, $\alpha=\frac{\pi}{2}$}
\end{figure}

\subsection{Isoptic curves of parabolas}
We study the isoptic curves of the generalized parabolas given by their equations in Theorem 3.3.

The algorithm described in the previous subsection can be repeated for further
conics as well. The difference is only in the equation
of the {\it compound isoptic} is because of the different conic equation.

It is clear that the coordinates of the tangent line $\Bsu=(u_1,u_2,1)$ and $\Bsv=(v_1,v_2,1)$ will be different from (8).

\begin{equation}
\begin{array}{l}
   u_1=-\frac{a x^2 (b+y-1)+y \sqrt{a b^2 x^2 \left(a x^2+b \left(y^2-1\right)+(y-1)^2\right)}}{a (b-1) x^3+b^2 x y^2}\\
   u_2=\frac{a x^2-b^2 y-\sqrt{a b^2 x^2 \left(a x^2+b \left(y^2-1\right)+(y-1)^2\right)}}{a (b-1) x^2+b^2 y^2}\\
   \\
   v_1=\frac{-a x^2 (b+y-1)+y \sqrt{a b^2 x^2 \left(a x^2+b \left(y^2-1\right)+(y-1)^2\right)}}{a (b-1) x^3+b^2 x y^2}\\
   v_2=\frac{a x^2 - b^2 y+ \sqrt{a b^2 x^2 \left(a x^2+b \left(y^2-1\right)+(y-1)^2\right)}}{a (b-1) x^2+b^2 y^2}\\
\end{array}
\end{equation}

\begin{theorem}
Let a parabola be given by its equation $ax^2+(b+1)y^2-2y=b-1$ (see Theorem 3.3). Then the compound $\alpha$-isoptic curve $(0<\alpha<\pi)$
of the considered conic has the equation
$$
\left(a \left(b \left(2 x^2+y^2-1\right)+(y-1)^2\right)+b^2 \left(y^2-1\right)\right)^2 \left| (y-1)^2 \left((y+1)^2 b^4-\right.\right.$$
$$
\left.\left.-2 a \left(2x^2+y^2+b (y+1)^2-1\right) b^2+a^2 \left((y-1)^2+b^2 (y+1)^2+2 b \left(2 x^2+y^2-1\right)\right)\right)\right|^{-1}=
$$
\newline
\[
=\left\{\begin{array}{r}
\begin{array}{ll}
 \cosh^2 (\alpha ), &a b^2 x^2 \left(a x^2+b \left(y^2-1\right)+(y-1)^2\right)\geq 0~\land \\
 & \left(1-u_1^2-u_2^2\right) \left(1-v_1^2-v_2^2\right)>0\land\\
 &x^2+y^2>1 \\
 & \\
 \cos^2 (\alpha ), & a b^2 x^2 \left(a x^2+b \left(y^2-1\right)+(y-1)^2\right)\geq 0~\land\\
 & x^2+y^2<1\\
 & \\
 \sinh^2 (\alpha ), & a b^2 x^2 \left(a x^2+b \left(y^2-1\right)+(y-1)^2\right)\geq 0~\land \\
 & \left(1-u_1^2-u_2^2\right) \left(1-v_1^2-v_2^2\right)<0,
\end{array}
\end{array}\right.
\]
wherein $u_{1,2}$ and $v_{1,2}$ are derived by \textnormal{(9)}.
\end{theorem}
The isoptic curves of the parabolas can be seen in Figure 4-6. The same convention has been used as on the previous figures.

\begin{figure}[p]
\centering
\includegraphics[scale=0.65]{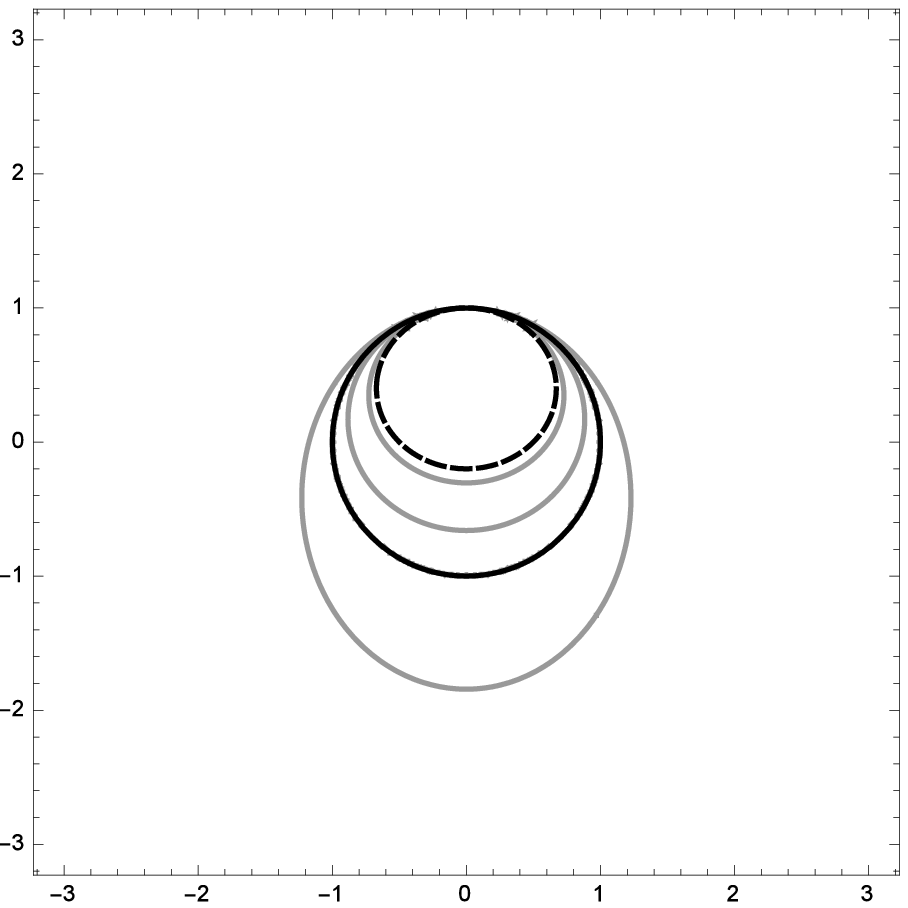}\includegraphics[scale=0.65]{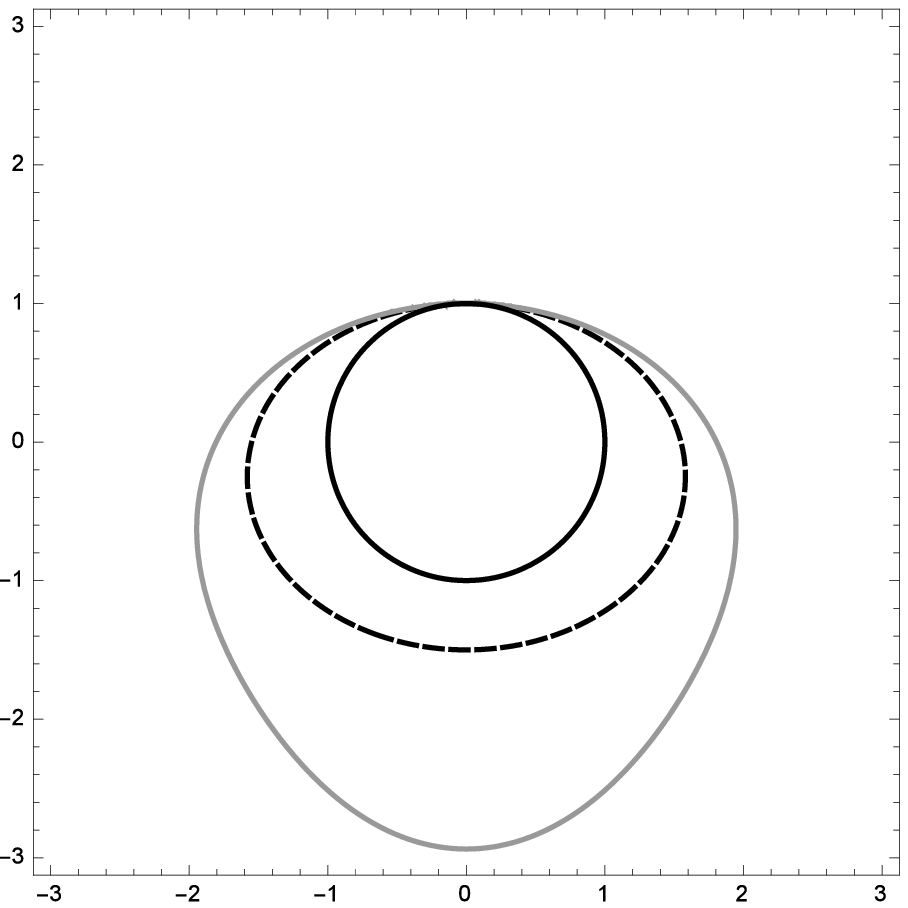}
\caption{Elliptic parabola(left): $a=2$, $b=1.5$, $\alpha=\frac{\pi}{3}$,\newline Parabola enclosing the absolute (right): $a=-2.5$, $b=-5$, $\alpha=\frac{7\pi}{18}$}
\includegraphics[scale=0.65]{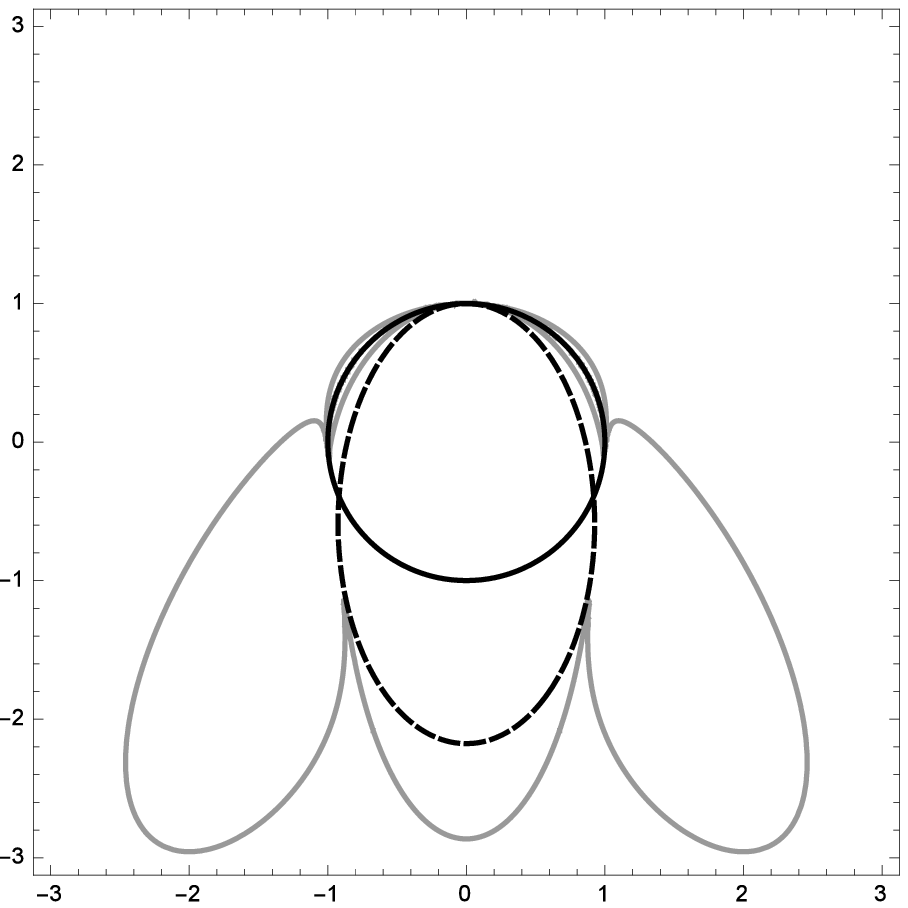}\includegraphics[scale=0.65]{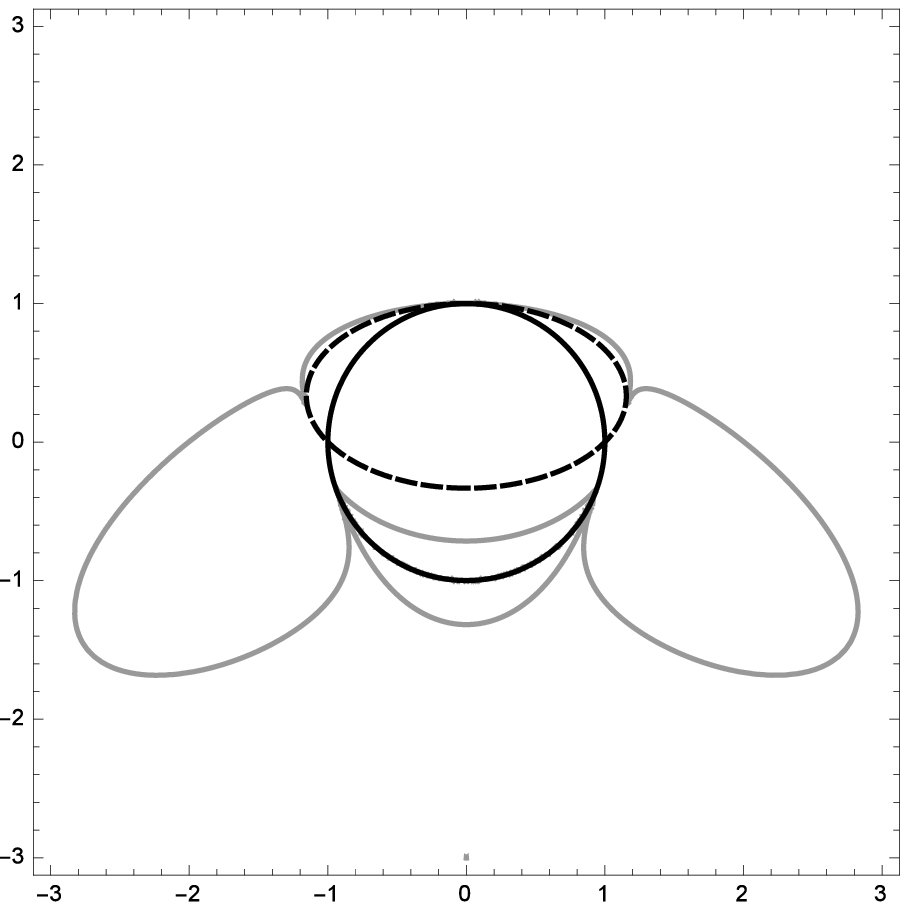}
\caption{Two sided parabola(left): $a=-5$, $b=-2.7$, $\alpha=\frac{\pi}{2}$, \newline Concave hyperbolic parabola(right): $a=1$, $b=2$, $\alpha=\frac{\pi}{2}$}
\includegraphics[scale=0.65]{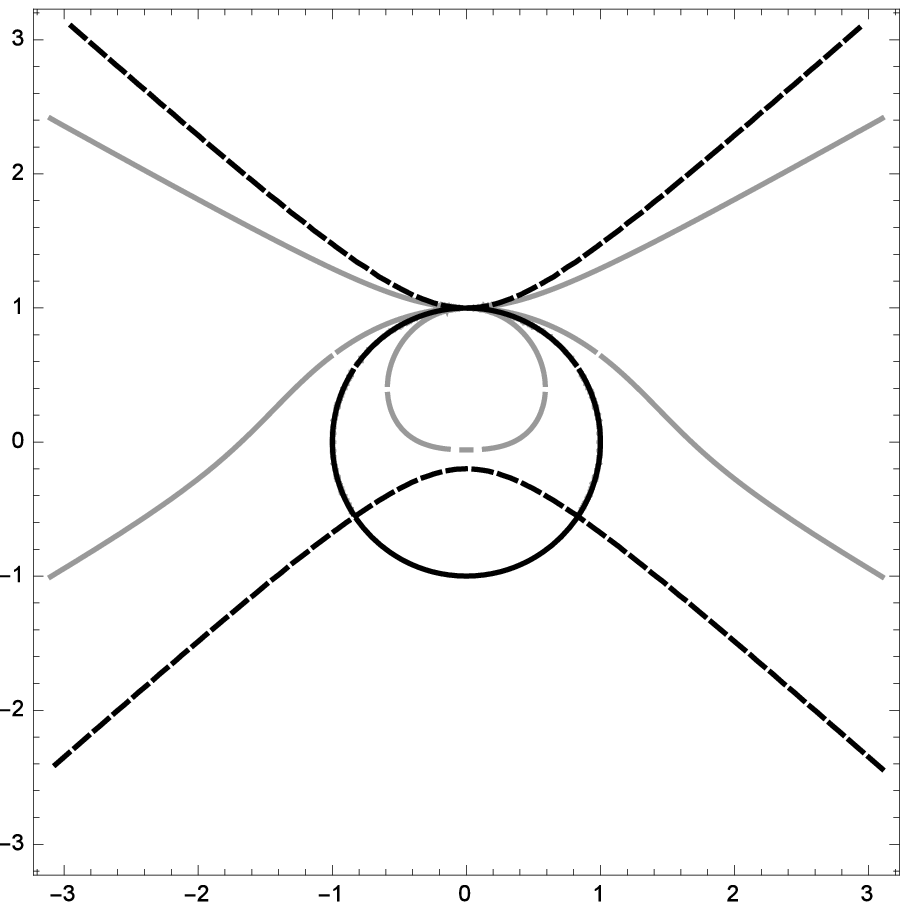}\includegraphics[scale=0.65]{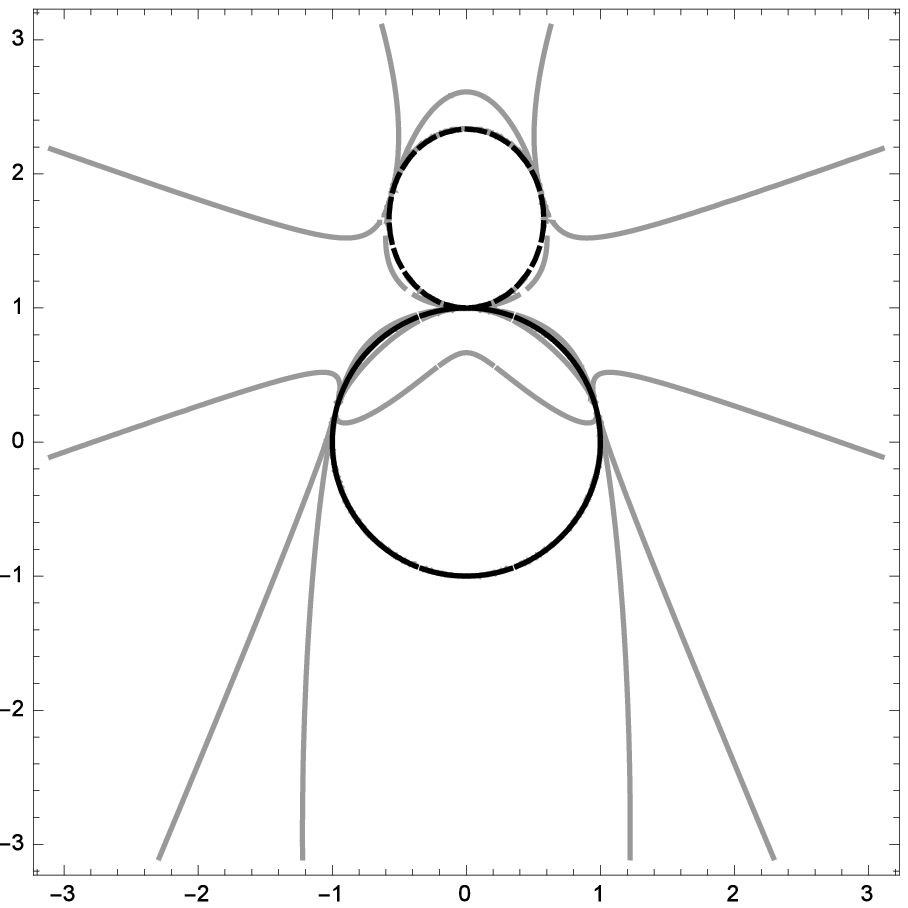}
\caption{Convex hyperbolic parabola(left): $a=-2$, $b=1.5$, $\alpha=\frac{\pi}{3}$,\newline Parabola excluding the absolute (right): $a=0.8$, $b=-0.4$, $\alpha=\frac{\pi}{3}$}
\end{figure}

\subsection{Isoptic curves of semi-hyperbola}
\begin{theorem}
Let the semi-hyperbola be given by its equation $a x^2+2 b y^2-2 y=0$, where $\left|b\right|<1
$ (see Theorem 3.2). Then the compound $\alpha$-isoptic curve
$(0<\alpha<\pi)$ of the considered conic has the equation
$$
\left(2 a \left(b \left(x^2+y^2\right)-y\right)+y^2-1\right)^2 \left| y^4+4 a^2 \left(x^2+y^2\right) \left(\left(b^2-1\right) x^2+(b
   y-1)^2\right)-\right.
$$
$$
\left.-4 a \left(y-\left(2 x^2+y^2\right) y+b \left(y^4+\left(x^2-1\right) y^2+x^2\right)\right)-2 y^2
+1\right|^{-1}=
$$
\newline
\[
=\left\{\begin{array}{r}
\begin{array}{ll}
 \cosh^2 (\alpha ), &a x^2 \left(a +2 y (b y-1)\right)\geq 0\land \\
 & \left(1-u_1^2-u_2^2\right) \left(1-v_1^2-v_2^2\right)>0\land\\
 &x^2+y^2>1 \\
 & \\
 \cos^2 (\alpha ), & a x^2 \left(a +2 y (b y-1)\right)\geq 0\land\\
 & x^2+y^2<1\\
 & \\
 \sinh^2 (\alpha ), & a x^2 \left(a +2 y (b y-1)\right)\geq 0\land \\
 & \left(1-u_1^2-u_2^2\right) \left(1-v_1^2-v_2^2\right)<0,
\end{array}
\end{array}\right.
\]
wherein
\begin{equation}
\begin{array}{l}
   u_1=-\frac{a x+\sqrt{a \left(a x^2+2 y (b y-1)\right)}}{y}\\
   u_2=\frac{a x^2-y+x\sqrt{a \left(a x^2+2 y (b y+1)\right)}}{y^2}\\
   \\
   v_1=\frac{-a x+\sqrt{a \left(a x^2+2 y (b y-1)\right)}}{y}\\
   v_2=\frac{a x^2-y-x\sqrt{a \left(a x^2+2 y (b y+1)\right)}}{y^2}.\\
\end{array}
\notag
\end{equation}
\end{theorem}
The isoptic curve of the semi-hyperbola can be seen in Figure 7.
\begin{figure}[ht]
\centering
\includegraphics[scale=0.65]{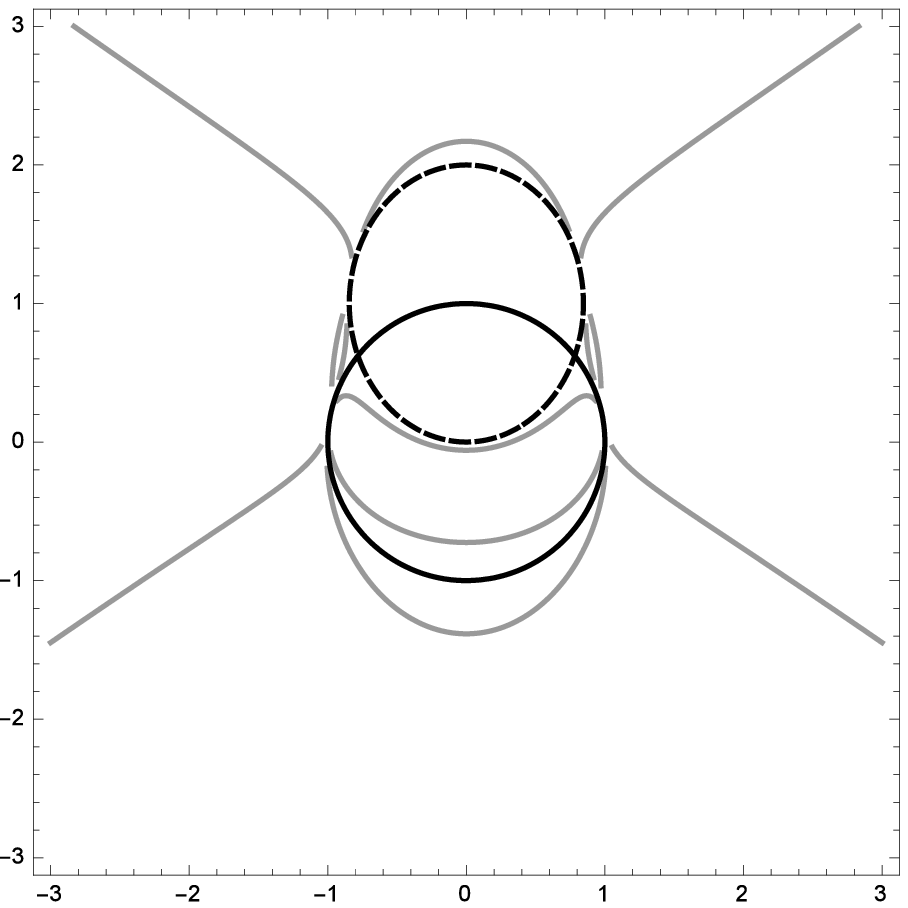}\includegraphics[scale=0.65]{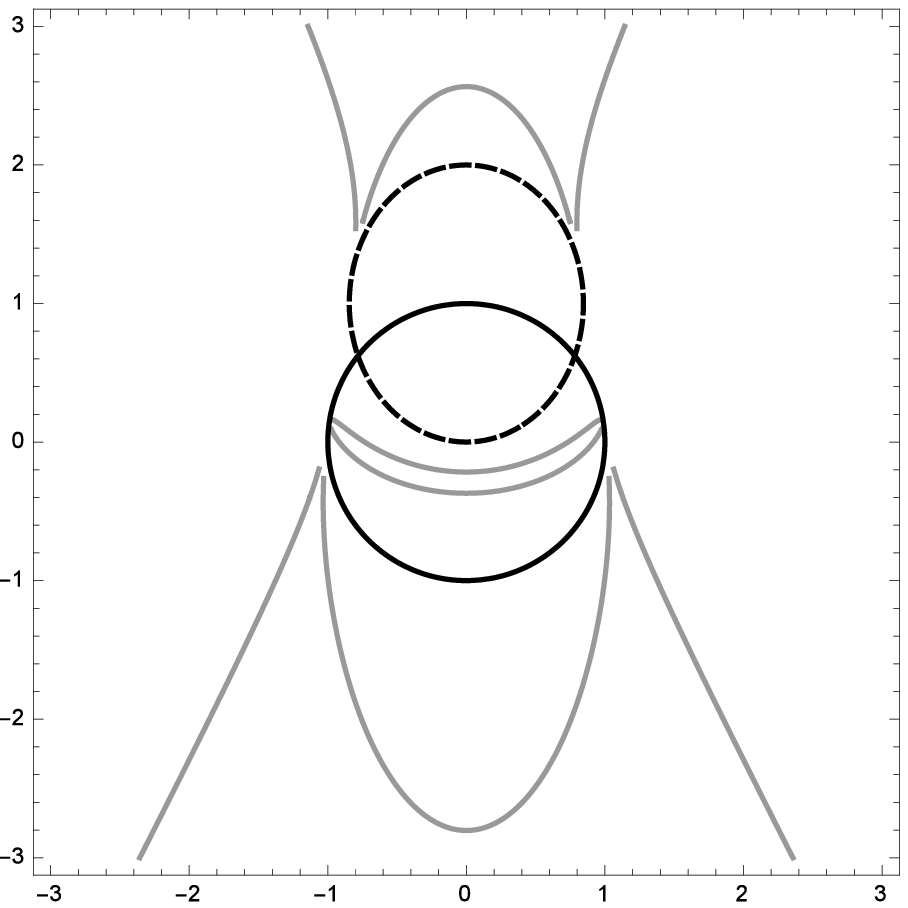}
\caption{Semi-hyperbola: $a=1.4$, $b=0.5$, $\alpha=\frac{\pi}{4}$ (left) and $\alpha=\frac{8\pi}{18}$ (right).}
\end{figure}

\subsection{Isoptic curves of the osculating parabola}
\begin{theorem}
Let the osculating parabola be given by its equation $\left(1-x^2-y^2\right) + 2 a (x+1) y=0$ (see Theorem 3.4).
Then the compound $\alpha$-isoptic curve $(0<\alpha<\pi)$ of the considered conic has the equation
$$
\frac{\left(-2 \left(x^2+y^2-1\right)+2 a (x+1) y+a^2 (x+1)^2\right)^2}{\left| a^2 (x+1)^3 \left(4(1-x)+4 a y +a^2
   (x+1)\right)\right| }=
$$
\newline
\[
=\left\{\begin{array}{r}
\begin{array}{ll}
 \cosh^2 (\alpha ), & \left(x^2+y^2-1-2 a (x+1) y\right)\geq 0\land \\
 & \left(1-u_1^2-u_2^2\right) \left(1-v_1^2-v_2^2\right)>0\land\\
 &x^2+y^2>1 \\
 & \\
 \cos^2 (\alpha ), & \left(x^2+y^2-1-2 a (x+1) y\right)\geq 0\land\\
 & x^2+y^2<1\\
 & \\
 \sinh^2 (\alpha ), & \left(x^2+y^2-1-2 a (x+1) y\right)\geq 0\land \\
 & \left(1-u_1^2-u_2^2\right) \left(1-v_1^2-v_2^2\right)<0,
\end{array}
\end{array}\right.
\]
wherein
\begin{equation}
\begin{array}{l}
   u_1=\frac{-(1+a y) (x-a y)+\sqrt{y^2 \left(x^2+y^2-1-2 a (x+1) y\right)}}{\left(x^2+y^2\right)-2 a x y+a^2 y^2}\\
   u_2=-\frac{ y^2-a  x (x+1) y+a^2 (x+1) y^2+x \sqrt{y^2 \left(x^2+y^2-1-2 a (x+1) y\right)}}{y \left(\left(x^2+y^2\right)-2 a x y+a^2 y^2\right)}\\
   \\
   v_1=-\frac{(1+a y) (x-a y)+\sqrt{y^2 \left(x^2+y^2-1-2 a (x+1) y\right)}}{\left(x^2+y^2\right)-2 a x y+a^2 y^2}\\
   v_2=-\frac{ y^2-a x (x+1) y+a^2 (x+1) y^2-x \sqrt{y^2 \left(x^2+y^2-1-2 a (x+1) y\right)}}{y \left(\left(x^2+y^2\right)-2 a x y+a^2 y^2\right)}.\\
\end{array}
\notag
\end{equation}
\end{theorem}
The Figure 8 shows some cases of the isoptic curve for the osculating parabola.
\begin{figure}[ht]
\centering
\includegraphics[scale=0.65]{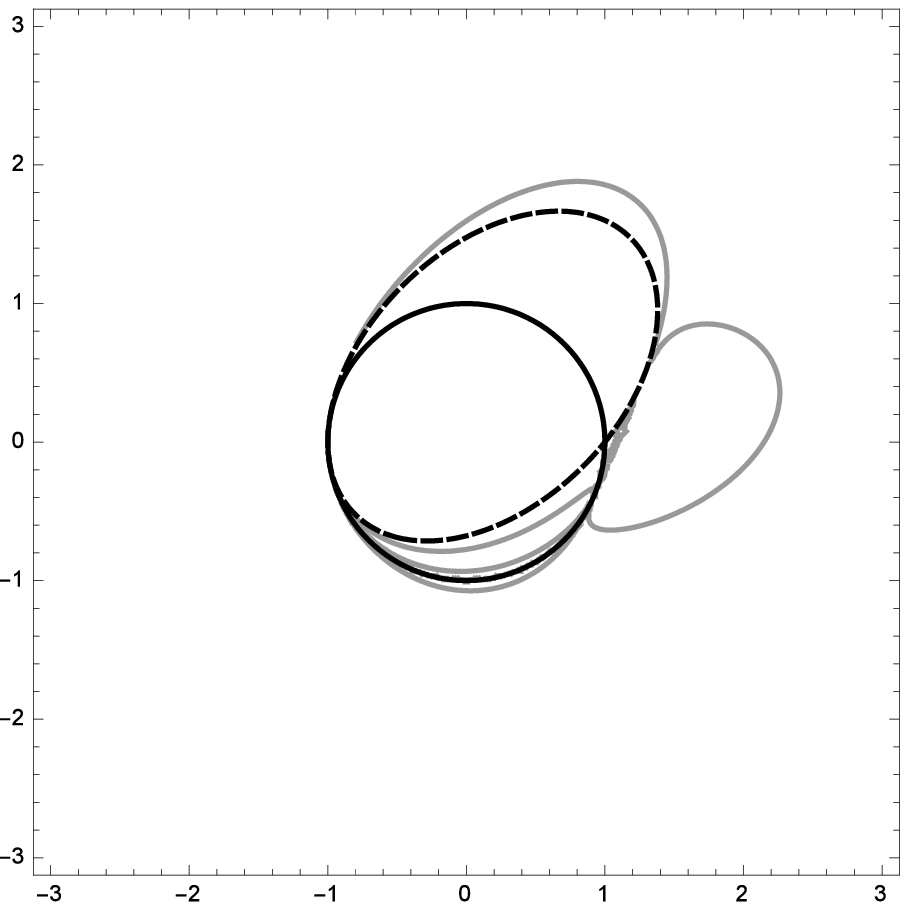}\includegraphics[scale=0.65]{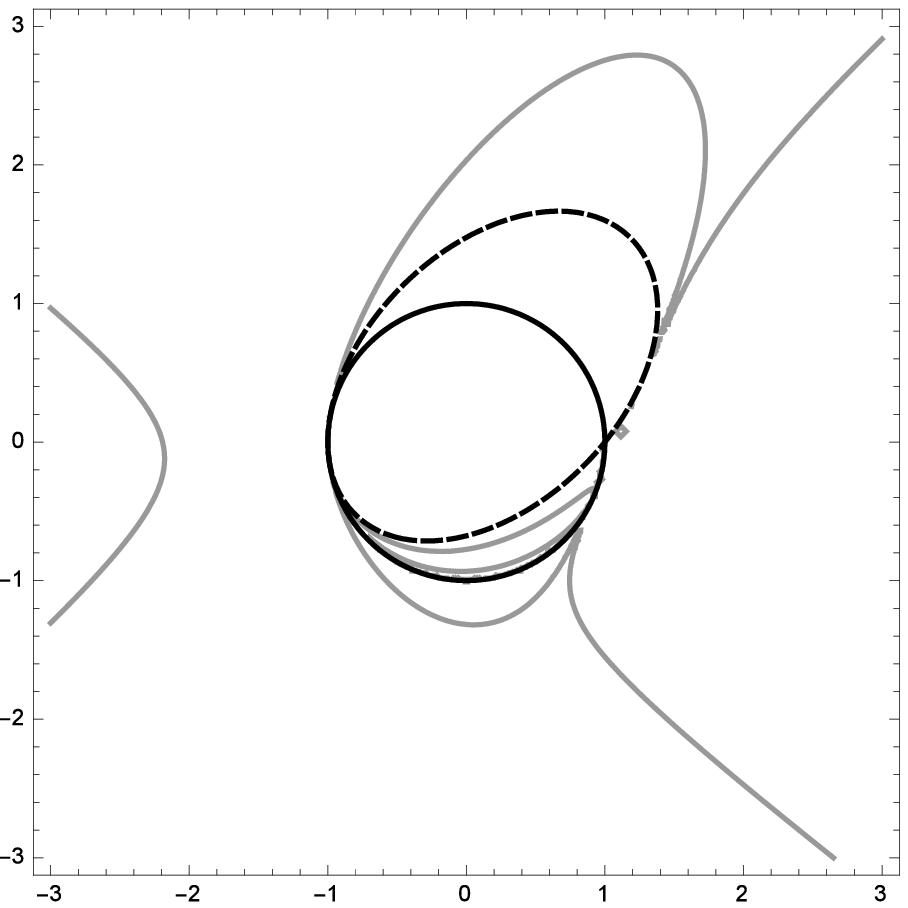}
\caption{Osculating parabola: $a=0.4$, $\alpha=\frac{\pi}{3}$ (left) and $\alpha=\frac{2\pi}{3}$ (right).}
\end{figure}

Our method is suited for determining the isoptic curves to generalized conic sections for all possible parameters. Moreover with this procedure above we may be able to determine the isoptics for other curves as well. Authors now turn attention toward 3D generalization. Already, some results in the Euclidean space can be seen in \cite{CsSz4}.

\section*{Acknowlidgement}
The authors would like to thank Professor Emil Moln\'ar for his very helpful discussions, instructions and comments to this paper, especially his constructive suggestions to the classification.


\end{document}